\documentclass[10pt,a4paper]{amsart}

\usepackage{pdfsync}
\usepackage{parskip}
\usepackage{amsmath}
\usepackage{setspace}
\RequirePackage{fix-cm}
\usepackage{biblatex} 
\usepackage{hyperref}
\usepackage[T1]{fontenc}
\usepackage[charter,cal=cmcal]{mathdesign}
\usepackage{xcolor}

\usepackage{mathptmx} 
\usepackage{mathrsfs}
\usepackage{mathalfa}
\usepackage[scaled=0.90]{helvet} 
\usepackage{courier} 
\normalfont
\usepackage[T1]{fontenc}
\usepackage{mwe} 
\usepackage{wrapfig}

\usepackage{marvosym}
\usepackage{graphicx}
\usepackage{epstopdf}
\usepackage[utf8]{inputenc}
\usepackage{lscape}
\usepackage{pifont}
\usepackage{tikz-cd}

\usepackage[curve,tips]{xypic}
\SelectTips{eu}{11}
\UseTips
\usepackage{hyperref}
\usepackage{comment}

\renewcommand{\epsilon}{\varepsilon}

\renewcommand{\emptyset}{\varnothing}

\newtheorem{theorem}{Theorem}[section]
\newtheorem{proposition}[theorem]{Proposition}
\newtheorem{corollary}[theorem]{Corollary}
\newtheorem{lemma}[theorem]{Lemma}

\theoremstyle{definition}
\newtheorem{example}[theorem]{Example}

\newtheorem{definition}[theorem]{Definition}

\newtheorem*{axiom}{Axiom}

\theoremstyle{remark}
\newtheorem{remark}[theorem]{Remark}

\usepackage{amsmath}

\title{Martin's axiom}
\author{Helena Jorquera Riera}

\begin{document}
\hfill\includegraphics[scale=0.2]{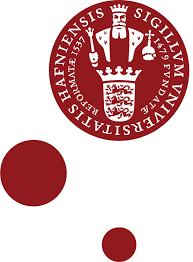}

\maketitle
\begin{center}
    
{\small 
Supervisor: Asger Dag Törnquist\\
University of Copenhagen\\
19$^{\textrm{th}}$ January 2023}
\end{center}
\begin{abstract}
    
     The axioms of ZFC provide a foundation for mathematics, however, there are statements independent of ZFC, such as the Continuum Hypothesis (CH). We discuss Martin's axiom, which is an alternative to CH that roughly states that if there is a cardinal strictly between $\omega$ and $2^{\omega}$ it "behaves" like $\omega$.
    
\end{abstract}
\section*{Introduction}

As a consequence of the (naive) set-theoretic paradoxes such as \textit{Russell's paradox} \cite{russells}, in the early 20th century Ernst Zermelo and Adolf Fraenkel formulated \textit{Zermelo–Fraenkel set theory} (ZF) \cite{zfc} as means of formalising set theory by using formal logic in place of natural languages in the theory thereby evading the paradoxes. ZF with the Axiom of Choice (ZFC) is a collection of axioms, i.e., statements which may be taken to be true without proofs, most of which are considered intuitive and generally accepted, although there is some controversy with the Axiom of Choice. From these axioms, all current mathematics results can be derived.

 George Cantor proved in 1891 \cite{Cantor} that there is more than one "size" of infinity by showing that "there are more" real numbers, $2^{\omega}$, than there are \textit{natural numbers}, $\omega$. He then tried to figure out whether there were any "sizes" between $\omega$ and $2^{\omega}$ but failed to do so. The \textit{Continuum Hypothesis} is the statement that there is no "size" between $\omega$ and $2^{\omega}$. This hypothesis remained unsolved for many years and was even put first on the list of Hilbert's problems in 1900. In 1963 Paul Cohen was awarded a Fields Medal for proving, with the help of Kurt Gödel, the hypothesis to be independent of ZFC \cite{cohen}.



If CH is false in some axiomatic system, then questions about the sizes $\kappa$ between $\omega$ and $2^{\omega}$ arise. Martin's Axiom (MA), introduced by Donald A. Martin and Robert M. Solovay \cite{martin}, is not as apparently intuitive as the other axioms of ZFC, but rather is a more ad hoc axiom in the sense that MA has the "good" consequences of the CH without the "bad" ones, hence we might want to add MA to ZFC as an alternative to CH rather than because of an intrinsic intuitive appearance. This is possible since MA is independent of ZFC, and since CH  is independent of MA, CH will remain unresolved if we accept MA, whereas CH trivially implies MA. Roughly Martin's Axiom says that the $\kappa$ between $\omega$ and $2^{\omega}$ "behave" like $\omega$, more explicitly, the answers to Questions 1-4 in Section \ref{sectionmartin} are "yes" for $\kappa=\omega$ and if we assume $MA(\kappa)$, then they are also "yes" for $\omega\leq \kappa<2^{\omega}$. 

We will give different equivalent forms of MA. The first one of which roughly states that if we have a collection of "not too many requirements", then we can find a \textit{filter} which ensures all requirements are satisfied. We will show MA to be equivalent to its restriction to \textit{Boolean Algebras}, and to a purely topological statement.

All results are taken from \cite{Kunen} unless stated otherwise.

\fontsize{11.3pt}{0.6cm}\selectfont

Note on notation: In what follows we use the symbols of first-order logic $\wedge, \exists, \forall, (,),\neg, =,$ and $v_i$, and others like $\in,<$. However, we often will use other symbols which may be written using combinations of the ones listed, e.g. $\Rightarrow$. We don't define the meaning of each symbol, we take it to be the "standard" one.

\section{The axioms of ZFC}
In this section, we give some set theory background; we list the ZFC axioms and define relations and orderings, which are necessary for the following sections.

\begin{definition}

    The \textit{universal closure} of a formula $\phi$ is the formula obtained from $\phi$ by universally quantifying all free variables appearing in $\phi$.
\end{definition}

We list the axioms of ZFC in the language of first-order logic:
\begin{axiom}
 Existence. $$\exists x(x=x).$$

 That is, there exists a set.
\end{axiom}
\begin{axiom}
 Extensionality. $$\forall x \forall y  (x=y \iff \forall z(z\in x \iff z\in y)).$$

 That is, a set is determined by the elements it contains.
\end{axiom}
\begin{axiom}
    Foundation. $$\forall x(\exists y(y\in x)\Rightarrow \exists y(\neg \exists z (z\in x \wedge z\in y))).$$

    That is, no set contains an infinite sequence of sets each containing the next. This avoids having self-referencing sets such as $x=\{x\}$.
\end{axiom}
\begin{axiom}
 Comprehension Scheme. The universal closure of the formula \\ $$\forall z(\exists y \forall x (x\in y \iff (x\in z \wedge \phi))),$$ where $y$ is not free in $\phi$.

That is, any subcollection of elements satisfying some property which doesn't refer to the collection itself is a set.

We write the unique $y$ satisfying the axiom as $\{x\in z: \phi\}$.
\end{axiom}
\begin{axiom}
 Pairing. $$\forall x \forall y \exists z(x\in z \wedge y\in z).$$
\end{axiom}
That is, we can make a set containing other sets that we already have. By Comprehension, we can define the set containing precisely $x$ and $y$, denoted $\{x,y\}$.
\begin{axiom}
 Union. $$\forall F\exists A\forall Y\forall x((x\in Y \wedge Y\in F) \Rightarrow x\in A).$$
 
 That is, given any family $F$ of sets $Y$ there is a set containing the elements of the elements of $F$, then by Comprehension, we can define the \textit{union} $\bigcup F$ of $F$ by $\{x\in A:\exists Y\in F(x\in Y)\}$ for some $A$ satisfying the axiom. 
\end{axiom}
\begin{axiom}
 Replacement scheme. The universal closure of $$\forall x\in A \exists ! y\phi(x,y)\Rightarrow \exists Y \forall x\in A \exists y\in Y \phi(x,y).$$%

 That is, if we associate one set $y$ to each member of a family of sets, we can make a set $Y$ out of all the $y$. The uniqueness of such a $y$ for each $x$ ensures that $Y$ is not "too big" to be a set.
\end{axiom}

\begin{definition}
    By Pairing, we can define the set $\{x\}$ whose only element is $x$. 
    
    The \textit{successor} $x$ is $S(x):=x \cup \{x\}$.

\end{definition}
\begin{definition}
    By Comprehension, we can define $0:=\{x\in y: x\neq x\}$ for any set $y$. Clearly $0=\emptyset$ is the set with no elements.
\end{definition}
\begin{axiom}
 Infinity. $$\exists x(0\in x\wedge \forall y\in x (S(y)\in x)).$$
 
 That is, there exists an infinite set.
\end{axiom}
\begin{axiom}
Power Set. $$\forall x\exists y \forall z(z\subseteq x\rightarrow z\in y).$$

That is, there is a set containing all subsets of a set. For any $y$ satisfying the Power Set Axiom the \textit{power set} of $x$ is $P(x)=\{z\in y:z\subseteq x\}$. This definition is justified by Comprehension.
\end{axiom}

\begin{axiom}
 Choice. $$\forall A\exists R( R \textrm{ well-orders } A),$$ where $R$ well-orders $A$ if $\langle A, R\rangle$ is a well-ordering in the sense of Definition \ref{defpartord}.

 This is equivalent to the statement that for any family of non-empty sets there is a set containing one element from each set.
\end{axiom}
Although not all axioms are necessary to all proofs, in particular, we never use the Axiom of Foundation, we assume the axioms of ZFC henceforth unless stated otherwise.


\begin{theorem}\label{thmsetofsets}
    There is no set of all sets.
\end{theorem}
\begin{proof}
    If the set $S$ of all sets existed, then by Comprehension, we could define $X=\{x\in S: x\notin S\}$. If $X\in S$, then $X \notin S$ and if $X\notin S$, then $X\in S$, a contradiction either way. This is Russell's paradox \cite{russells}.
\end{proof}

\begin{definition}

    The \textit{intersection} $\cap F$ of the sets in a family $F$ is defined to be $\{x\in X:\forall Y\in F(x\in Y)\} $ for $F\neq 0$ and some (any) $X\in F$.

    Here we are restricting $x$ to be in some already existing set $X\in F$ so that we can apply Comprehension. Hence we need $F\neq 0$. Indeed if $F=0$, then every set $x$ satisfies $\forall Y\in F(x\in Y)$, but no set of all sets exists by Theorem \ref{thmsetofsets}. 
\end{definition}

By using Pairing twice, we can define the following.
\begin{definition}
         The \textit{ordered pair} of $x$ and $y$ is $\langle x,y\rangle= \{x,\{x,y\}\}$.
\end{definition}

\begin{definition}
    The Cartesian product of $A$ and $B$ is $A\times B=\{\langle a,b\rangle: a\in A \wedge b\in B$\}.
\end{definition}
For each $b\in B$ we can define $X_b:=\{\langle a,b\rangle: a\in A\}$ by Replacement. Then we can define $\{X_b:b\in B\}$ again by Replacement. Now by Union, we can define $\bigcup\{X_b:b\in B\}$, which is equal to $A\times B$.

\begin{definition}
    A \textit{relation} $R$ is a collection of ordered pairs. If $R$ is a relation, then the \textit{domain} of $R$ is $Dom(R):=\{x:\exists y(\langle x,y \rangle \in R)\}$ and its \textit{range} is $Ran(R):=\{x:\exists y(\langle y,x \rangle \in R)\}$.

    A function $f$ is a relation s.t. if $\langle x,y \rangle \in f $ and $\langle x,y' \rangle \in f$, then $y=y'$.

    If $dom(f)=A$, $ran(f)=B$ and $C\subseteq A$, then $f|_C=f\cap C\times B$ is the restriction of $f$ to $C$, and $f''C=ran(f|_C)$.
\end{definition}

We usually think of functions as associating a single value in its range to each value in its domain; if $\langle x,y \rangle \in f$ we say that $y$ is the value of $f$ at $x$.
\begin{definition}\label{defpartord}
    A \textit{partial order} $\langle P, \leq\rangle$ (or just P) is a set $P\neq 0$ together with a transitive and reflexive relation $\leq$, i.e. $$\forall x,y,z\in P (x\leq y \wedge y\leq z \Rightarrow x\leq z)$$ and $$\forall x\in P(x\leq x).$$   

   A \textit{total ordering} $\langle A, R\rangle$ is a pair consisting of a set $A$ and a relation $R$ such that transitivity: $$\forall x,y,z\in A (xR y \wedge yR z \Rightarrow xR z),$$ trichotomy: $$\forall x,y\in A(x=y\vee xRy \vee yRx)$$ and irreflexivity: $$\forall x\in A (\neg xRx)$$ hold.

    A \textit{well-ordering} is a total ordering $\langle A, R\rangle$ s.t. every non-0 subset of $A$ has an $R$-least element.
\end{definition}
\begin{example}\label{exfam}
    The set of the entire world's population with ordering given by $x\leq y$ if $x$ is a descendant of $y$ is a partial order $\langle P , \leq \rangle$.
\end{example}
\begin{example}
    Whereas, the partial order $(\mathbb{R},\leq)$ with set the real numbers and where $x\leq y$ if $\exists n>0\in \mathbb{N}(x+n=y)$ is not a partial order since trichotomy fails.

\section{Ordinals}
In this section, we introduce the \textit{ordinals} and show that any well-ordering is \textit{isomorphic} to an \textit{ordinal}. Then we define some operations on them.

\begin{definition}
   A set $x$ is \textit{transitive} if every element of $x$ is a subset of $x$.
\end{definition}
\begin{example}
    The sets $\{\{0\},0\}$ and, if we don't assume Foundation, $x=\{0,\{0\},x\}$ are transitive sets.  
\end{example}
\begin{definition}
    An ordinal $x$ is a transitive set s.t. $\langle x, \in_x\rangle$ is a well-ordering, where $\in_x=\{\langle y,z\rangle \in x\times x: y\in z\}$.

    Throughout we use greek letters to denote ordinals except for $\phi$ and $\psi$ which denote functions or formulas (which one is being used will be clear by context).
\end{definition}
\begin{example}
   The set $\{\{0\},0\}$ is an ordinal but even if we don't assume Foundation, $x=\{0,\{0\},x\}$ isn't.
\end{example}
\begin{theorem}\label{thmordoford}
    There is no set of all ordinals.
\end{theorem}
\begin{proof}
    If the set of all ordinals $ON$ were to exist it would be an ordinal and a member of itself, contradicting the irreflexivity of $\in$. 
\end{proof}
\begin{definition}
    For any sets $A,B$ and relations $R,S\subseteq A\times A, B\times B$ respectively, we say $\langle A,R\rangle$ and $ \langle B,S \rangle$ are \textit{isomorphic} and denote it $\langle A,R\rangle\cong \langle B,S \rangle$ if there exists a bijection $f:A\rightarrow B$ with $\forall a,b \in A(aRb\iff f(a)Sf(b))$ called an \textit{isomorphism}.
\end{definition}
\begin{definition}
    Let $\langle A, R \rangle $ be a well-ordering, we define $pred(A,x,R):= \{a\in A: aRx\}$.
\end{definition}
\begin{lemma}\label{pftransind}
    For any well-ordering $\langle A, R \rangle $ we have $\langle A,R \rangle \not\cong \langle pred(A,x,R), R\rangle$.
\end{lemma}
This is my attempt at proving the previous lemma.
\begin{proof}

    Suppose $f: A\rightarrow pred(A,x,R)$ is an isomorphism.
    We show that if $\forall x<y (x=f(x)\Rightarrow f(y)=y)$ hence $\forall y\in A (f(y)=y)$ (this is a proof by transfinite induction, see Section \ref{secttrind}), a contradiction since $A\neq pred(A,x,R)$.
    
    Suppose there exists $y\in A$ s.t. $f(y)\neq y$ and let $y$ be the smallest with this property. Then every $z\in pred(A,x,R)$ s.t. $xRy$ is the image of $z\in A$, hence since $f$ is injective and $f(y)\neq y$ we have $yRf(y)$. 
    
    Also, $f$ is surjective so $y\in pred(A,x,R)$ is the image of some element $f^{-1}(y)\neq y$ in $A$. Since $f$ is a function we have $f^{-1}(y)\neq x$ for any $xRy$ thus $yRf^{-1}(y)$. 

Altogether, $y R f^{-1}(y)$ and $yRf(y)$, so $f$ doesn't preserve the ordering $R$. Hence, no such $y$ exists, thus $f(y)=y$.
\end{proof}

\begin{lemma}\label{pftransind2}
    If $\langle A,R \rangle\cong \langle B, S \rangle$ (well-orderings), then there is a unique isomorphism between them.
\end{lemma}
This is my attempt at proving the previous Lemma.

\begin{proof}
    Suppose $f,g: A\rightarrow B$ are two distinct such isomorphisms. Then there is smallest $y$ s.t. $w:=f(y)\neq g(y)=:z$. Suppose w.l.o.g. that $zSw$. Then since $f$ is surjective $z$ has a preimage $f^{-1}(z)$ in $A$. Since $f$ is a function $f^{-1}(z)\neq x$ for $xRy$ and for $x=y$. Hence $yRf^{-1}(z)$ and $zSf(y)$.
\end{proof}
\begin{theorem}\label{thm3options}
    For any two well-orderings $\langle A,R \rangle$ and $\langle B, S \rangle$ one and only one of the following holds:

    (1) $\langle A,R \rangle\cong \langle B, S \rangle$.\\
    (2) $\exists x\in B \langle A,R \rangle\cong \langle pred(B,x,S), S \rangle$.\\
    (3) $\exists y\in A\langle pred(A,y,R),R \rangle\cong \langle B, S \rangle$.
    
\end{theorem}

\begin{theorem}\label{thmCiso}
    If $\langle A, R\rangle$ is a well-ordering there is a unique ordinal $C$ s.t. $\langle A, R\rangle \cong \langle C,\in \rangle =: C$.
\end{theorem}
\begin{proof}
Uniqueness follows from the fact that if two ordinals $x$ and $ y$ are isomorphic, then $x=y$.

Let $B:=\{a\in A: \exists x(x \textrm{ is an ordinal }\wedge pred(A,a,R)\cong x)\}$ and define $f:B\rightarrow ran(f)=:C$ by $f(b)=\textrm{ the unique ordinal $x$ s.t. } pred(A,b,R)\cong x$. Then $f$ is an isomorphism, that is, $\langle B, R\rangle \cong C$ and $C$ is an ordinal. 

If $B\neq A$, then by Theorem \ref{thm3options}, there is some $a\in A$ s.t. $B=pred(A,a,R)$ but then $pred(A,a,R)\cong C$ so $b\in B$, contradicting that $R$ is irreflexive. Hence $B=A$ so $\langle A, R\rangle \cong C$. 
\end{proof}
Hence, we might want to restrict our study of well-orderings to ordinals:
\begin{definition}
   For $\langle A, R\rangle$ a well-ordering $\textrm{type}( A, R)$ is the unique ordinal $C$ s.t. $\langle A, R\rangle \cong  C$.  
\end{definition}

\begin{definition}
    For $X$ a set of ordinals $\textrm{sup}(X)=\bigcup X$ and if $X\neq 0$, $inf(X)=\bigcap X$.
\end{definition}


\begin{definition}

   An ordinal $\alpha$ is a \textit{successor} if $\exists \beta(S(\beta)=\alpha)$. If $\alpha\neq 0$ and is not a successor it is a \textit{limit ordinal}.

   We denote $1:=S(0)=\{0\}, 2:=S(1)=\{0,1\}=\{0,\{0\}\}, 3:=S(2)=\{0,1,2\}=\{0,\{0\},\{0,\{0\}\}\}$, etc.
\end{definition}

\begin{definition}
      An ordinal $\alpha$ is a \textit{natural number} if it is 0 or if every non zero $\beta\leq \alpha$ is a successor. 
    
    The set of all natural numbers $\omega$ exists by the Axioms of Infinity and Comprehension.
\end{definition}
\begin{proposition}
    The set $\omega$ is a limit ordinal.
\end{proposition}
\begin{proof}
    By definition, every $\beta <\omega$ is a successor or 0. Then if $\omega$ were a successor $\omega$ would be a natural number, hence $\omega\in \omega$, contradicting that $\in$ is a strict ordering.
\end{proof}
\begin{example}
    The successor $S(\omega)$ is not a natural number since $\omega\leq S(\omega)$ but $\omega$ is not a successor.

    More generally, the natural numbers are all the ordinals below $\omega$.
\end{example}

\begin{definition}
    We define the operation $+$ on ordinals by  $\alpha+\beta=\textrm{type}(\alpha\times\{0\}\cup \beta\times \{1\}, R)$ where $R=\{\langle\langle\xi, 0\rangle,\langle\eta,0\rangle\rangle:\xi<\eta<\alpha\}\cup\{\langle\langle\xi, 1\rangle,\langle\eta, 1\rangle\rangle:\xi<\eta<\beta\}\cup ((\alpha\times \{0\})\times(\beta\times \{1\})$.

    When some operation is only defined for e.g. ordinals we define it to be 0 when one of the elements in the operation is not an ordinal. 
\end{definition}
This definition is justified since we can define $\alpha +\beta=\{\gamma: \textrm{type}(\alpha\times\{0\}\cup \beta\times \{1\}, R)=\gamma\}$ this is justified by Replacement and Comprehension since by Theorem \ref{thmCiso} there is a unique such $\gamma$ for each pair $\alpha,\beta$.

\begin{example}
$1+2=\textrm{type}(1\times\{0\} \cup 2\times\{1\}, R)=\textrm{type}(\langle 0,0\rangle,\langle 0,1\rangle,\langle 1,1\rangle)$ with \\$R=\{\langle\langle0,0\rangle,\langle0,1\rangle\rangle,\langle\langle0,1\rangle,\langle\langle1,1\rangle\rangle, 
\langle\langle0,0\rangle,\langle1,1\rangle\rangle\}$.

Thus \{$\langle 0,0\rangle,\langle 0,1\rangle,\langle 1,1\rangle\}\cong\{0,1,2\}\cong 3$ so $1+2=3$.

\end{example}
\begin{lemma}\label{lem+} For any $\alpha$ and $\beta$ we have:\\
    (1) $\alpha+0=\alpha$.\\
    (2) $\alpha+S(\beta)=S(\alpha+\beta)$.\\
    (3) If $\beta$ is a limit ordinal $\alpha+\beta=\textrm{sup}\{\alpha+\gamma: \gamma<\beta\}.$
\end{lemma}
\begin{proof}[Proof of (1)]
      $\alpha+0=\textrm{type}(\alpha\times \{0\}\cup 0\times\{1\})=\textrm{type}(\alpha\times \{0\})=\alpha$.
\end{proof}
\begin{lemma}
    The operation + is not commutative.
\end{lemma}
This is my attempt at proving the previous lemma.
\begin{proof}
    We show $1+\omega=\omega\neq \omega+1$:

    $1+\omega=\textrm{type}(1\times\{0\}\cup \omega\times\{1\}, R)$ where ordering $\{1\times\{0\}\cup \omega\times\{1\}$ by $R$ gives \\$\{\langle 1,0\rangle, \langle 0, 1 \rangle, \langle 1, 1 \rangle, \langle 2, 1 \rangle, ...\}$. Hence there is an isomorphism $f:\langle 1\times\{0\}\cup \omega\times\{1\}, R\rangle\rightarrow \omega$ where $f(1,0)=0$ and $f(n,1)=n+1$ for all $n\in \omega$. Thus $1+\omega=\omega$.

    $\omega+1=\textrm{type}(\omega\times\{0\}\cup 1\times\{1\}, R)$ where ordering $\{\omega\times\{0\}\cup 1\times\{1\}$ by $R$ gives\\ $\{\langle 0,0\rangle, \langle 1,0 \rangle, \langle 2,0 \rangle,..., \langle 1, 1 \rangle\}$. Hence there is an isomorphism $f:\langle \omega\times\{0\}\cup 1\times\{1\}, R\rangle\rightarrow S(\omega)$ where $f(n,0)=n$ for all $n\in \omega$ and $f(1,1)=\omega$. Thus $\omega+1=S(\omega)$.
\end{proof}

\begin{definition}
    We define the operation $\cdot$ on ordinals by 
 $\alpha\cdot\beta=\textrm{type}(\beta\times \alpha, R)$ where $R=\{\langle \alpha, \beta \rangle, \langle \alpha', \beta' \rangle:  \alpha < \alpha' \vee (\alpha =\alpha'\wedge \beta< \beta')  \}$. 
\end{definition}
Notice that we change the order of $\alpha,\beta$ "inside" the type operation.
\begin{lemma}
    The operation $\cdot$ is non-commutative.
\end{lemma}
This is my attempt at proving the previous lemma.

\begin{proof}
    We show $2\cdot \omega=\omega\neq \omega\cdot 2$:

    $\omega\cdot 2 =\textrm{type}(2\times \omega, R)=\textrm{type}(\{\langle 0, 0\rangle, \langle 0, 1\rangle, \langle 0, 2\rangle, ..., \langle 1, 0\rangle,\langle 1, 1\rangle, \langle 1, 2\rangle, ...\}, R)$. When ordering $2\times \omega$ by $R$ we get $\{\langle 0, 0\rangle, \langle 0, 1\rangle, \langle 0, 2\rangle, ..., \langle 1, 0\rangle,\langle 1, 1\rangle, \langle 1, 2\rangle, ...\}\cong \omega +\omega \not\cong \omega$. Thus $\omega\cdot 2\neq \omega$.

$2\cdot \omega =\textrm{type}(\omega\times 2, R)=\textrm{type}(\{\langle 0, 0\rangle, \langle 1,0\rangle, \langle 2,0\rangle, ..., \langle 0,1\rangle,\langle 1,1\rangle, \langle 2,1\rangle, ...\}, R)$ where ordering $2\times \omega$ by $R$ gives $\{\langle 0, 0\rangle, \langle 0, 1\rangle, \langle 1,0\rangle, \langle 1,1\rangle,\langle 2,0\rangle, \langle 2,1\rangle, ...\}\cong \omega$. Thus $2\cdot \omega=\omega$ %
\end{proof}

\begin{lemma}
    For any $\alpha,\beta, \gamma$ we have:\\
    (1) $\alpha\cdot(\beta\cdot\gamma)=(\alpha\cdot\beta)\cdot\gamma $.\\
    (2) $\alpha\cdot 0=0$.\\
    (3) $\alpha \cdot 1=\alpha$.\\
    (4) $\alpha\cdot S(\beta)=\alpha\cdot\beta+\alpha$.\\
    (5) If $\beta$ is a limit ordinal $\alpha\cdot\beta=\textrm{sup}\{\alpha\cdot \xi: \xi<\beta\}$.\\
    (6) $\alpha\cdot(\beta+\gamma)=\alpha\cdot\beta+\alpha\cdot\gamma$.
\begin{remark}
    Distributivity on the right fails, e.g. $(1+1)\cdot \omega=2\cdot \omega =\omega \neq \omega+\omega$.
\end{remark}

\end{lemma}

\section{Transfinite Recursion}\label{secttrind}
In this section, we introduce \textit{Transfinite Recursion} which allows us to define functions on a given value given that we understand the function on the smaller values, which in turn allows us to define the function as a whole. An example of this is the definition of factorials; $n!=(n-1)!\cdot n$, that is, we know the value of the function sending a natural number to its factorial on $n$ if we know its value on $n-1$.

However, we want functions to take values on "collections" which need not constitute sets, e.g. we might want to define a function taking as values any ordinal and we know by Theorem~\ref{thmordoford}, that there is no set of all ordinals. Hence we introduce a more general notion of \textit{classes}, informally, a \textit{class} is a collection of the form $\{x:\phi\}$, with no restriction on $\phi$. Not all classes exist formally in ZFC but we use them to simplify the notation.
\begin{definition}
    $V=\{x:x=x\}$ and $ON=\{x: x\text{ is an ordinal}\}$.
\end{definition}

These classes are not sets by Theorems \ref{thmsetofsets} and \ref{thmordoford}, hence, statements like $x\in ON$ are not literal but rather must be read as "$x$ is an ordinal"

In the following, we will use proofs by \textit{transfinite induction}, which is a proof of a statement of the form $\forall x \psi(x)$ by showing that for any $\alpha$ we have $(\forall \beta<\alpha (\psi(\beta)))\Rightarrow \psi(\alpha)$. Then if $\exists \alpha (\neg \psi(\alpha))$, then there is a least such $\alpha$.  But then for every $\beta < \alpha$ we have $\psi(\beta)$ hence $\psi(\alpha)$, a contradiction.

\begin{example}
    The proofs of Lemmas \ref{pftransind} and \ref{pftransind2} are proofs by transfinite induction.
\end{example}


A function on $a$ can be defined by recursion from the information below $a$:
\begin{theorem}[Transfinite Recursion] If $F: V\rightarrow V$ there is a unique $G:ON\rightarrow V$ s.t. $\forall\alpha(G(\alpha)=F(G|_{\alpha}))$.
\end{theorem}

\begin{proof}
    Uniqueness: Assume that $\forall\alpha (G_1(\alpha)=F(G_1|_{\alpha}))$ and $\forall\alpha (G_2(\alpha)=F(G_2|_{\alpha}))$. Assume further that $\forall {\beta}<\alpha(G1({\beta})=G_2({\beta}))$ for a proof by transfinite induction,\\ $\Rightarrow \forall \beta<\alpha(F(G_1|_{\beta})=F(G_2|_{\beta}))$ \\ $ \Rightarrow \forall {\beta}<\alpha (\forall \gamma<{\beta}(F(G_1(\gamma))=F(G_2(\gamma))))$ \\ $\Rightarrow \forall {\beta}<\alpha (G_1({\beta})=G_2({\beta}))$ \\$ \Rightarrow \forall {\beta}<\alpha (F(G_1({\beta}))=F(G_2({\beta})))$\\ $\Rightarrow F(G_1|_{\alpha})=F(G_2|_{\alpha}) $\\ $\Rightarrow G_1({\alpha})=G_2({\alpha})$ \\$ \Rightarrow \forall\alpha (G_1({\alpha})=G_2({\alpha}))$ by Transfinite Recursion.

Existence: We say $g$ is a $\delta$-approximation if $\forall\alpha<\delta (g(\alpha)=F(g|_{\alpha}))$. Then if $g$ and $g'$ are $\delta$- and $\delta'$-approximations, then $g|_{(\delta \cap \delta)}=g'|_{(\delta \cap \delta)}$ (by the same argument as in the uniqueness proof). Suppose for all $\delta<\delta'$ there is a $\delta$-approx $g$. Then define $g'|_{\delta'}=g|_{\delta'}=g$ and $g'(\delta')=F(g)$ (which is well defined since approximations are unique). Then $g'$ is a $\delta'$-approx. Then we can define $G(\alpha)=g(\alpha)$ for $g$ the $\delta$-approximation for any $\delta>\alpha$.
\end{proof}
\begin{example}\label{exn0}
    For instance, we can use Transfinite Recursion to define the following: 


$n_0=0$\\
$n_{i+1}=n_{i}+i$ for all $i\in \omega$.

That is, $\{n_0,n_1,...\}=\{0,1,3,6,10, ...\}$.
\end{example}
We could also have defined $\alpha+{\beta}$ using recursion by taking Lemma~\ref{lem+} as a definition.

For another example, we use recursion to define ordinal exponentiation, which is distinct from cardinal exponentiation given in Definition \ref{defcardexp}:

\begin{example}\label{defordexp}(Ordinal exponentiation)\\
$\alpha^0=1$\\
$\alpha^{\beta+1}=\alpha^{\beta}\cdot \alpha $\\
If $\beta$ is a limit, then $\alpha^{\beta}=sup\{\alpha^{\xi}: \xi<\beta\}$.

\end{example}


\section{Cardinals}
In this section, we introduce \textit{cardinality}, which roughly tells us when two sets have the same "size". Then we introduce some results useful for computation, and we prove Cantor's Theorem \ref{thmcantor} that a set and its power set have different \textit{cardinalities}.

\begin{definition}
    We write $A\lesssim B$ if there is a 1-1 function from $A$ to $B$, and 
    we write  $A\approx B$ if there is a bijection from $A$ to $B$.
\end{definition}
\begin{theorem}
    $A\lesssim B, B\lesssim A \Rightarrow A\approx B$.
\end{theorem}

By Theorem \ref{thmCiso}, there is a unique ordinal isomorphic to any well-ordering, however, in general, there isn't a unique ordinal in bijection with a well-ordering, e.g., there is a bijection $\omega\rightarrow \omega+1$ given by $n\rightarrow n-1$ for every $n\neq 0$ and $0\rightarrow \omega$. This bijection is not an isomorphism, so $\omega$ and $\omega+1$ are different ordinals but since there is a bijection between them, we take them to have the same \textit{size}.  

\begin{definition}
    The \textit{cardinality} (or \textit{size}) $|A|$ of a well-ordered set $A$ is the least ordinal $\alpha$ s.t. $\alpha\approx A$.

    An ordinal $\alpha$ is a \textit{cardinal} if $\alpha=|\alpha|$. Henceforth we use $\kappa$ to denote a cardinal.

    A set $A$ is \textit{finite} if $|A|<\omega$ (otherwise it is \textit{infinite}) and it is \textit{countable} if $|A|\leq\omega$ (otherwise it is \textit{uncountable}).
\end{definition}
\begin{example}
    The ordinal $\omega$ is a cardinal since it is the smallest infinite ordinal hence it can't be in bijection with any smaller ordinal.
\end{example}
\begin{example}\label{exfinsubsets}
    Let $S$ be the set of all finite subsets of $\omega$, which exists by Power Set axiom and Comprehension, then $|S|={\omega}$; we order $S$ like so $S=\{\{0\},\{1\},\{0,1\},\{2\},\{0,2\},\{1,2\},\{0,1,2\},...\}$. Then we define a bijection $S\rightarrow \omega$ by $s\rightarrow n$ where $n$ is the position of $s\in S$ in our ordering.
%
\end{example}
\begin{theorem}\label{thmmaxkl}
    For any cardinals $\kappa$ and $\lambda$ we have $|\kappa+\lambda|=|\lambda+\kappa|=|\kappa\cdot\lambda|=|\lambda\cdot\kappa|=|\kappa\times\lambda|=\textrm{max}\{\kappa,\lambda\}$.
\end{theorem}
\begin{proof}
    We first show $|\kappa\times\kappa|=\kappa$ by transfinite induction. Assume $|\alpha\times\alpha|<\kappa$ for all $\alpha<\kappa$. Define an ordering by $\langle \alpha,\beta\rangle\ll \langle\gamma,\lambda\rangle$ iff 
    $$\textrm{max}\{\alpha,\beta\}<\textrm{max}\{\gamma,\lambda\}\vee(\textrm{max}\{\alpha,\beta\}=\textrm{max}\{\gamma,\lambda\}\wedge (\alpha<\gamma\vee(\alpha=\gamma \wedge \beta<\lambda)).$$

Then for each $\langle \alpha, \beta \rangle$ we have $|\{\langle \alpha',\beta'\rangle:\langle \alpha',\beta'\rangle\ll \langle \alpha,\beta\rangle\}|\leq|(\max\{\alpha,\beta\}+1)\times (\max\{\alpha,\beta\}+1)|<\kappa$ since $\kappa$ is a cardinal hence $\max\{\alpha,\beta\}+1<\kappa$. That is, each element has fewer than $\kappa$ smaller elements hence $\textrm{type}(\kappa \times \kappa, \ll)\leq \kappa$ thus $|\kappa\times\kappa|\leq \kappa$ and clearly $|\kappa\times\kappa|\geq \kappa$.

    Let $\kappa=max\{\kappa,\lambda\}$, then $|\kappa|=|\kappa\times\kappa|\geq |\kappa\times\lambda|$ and $|\kappa|\leq|\kappa\times\lambda|$ hence $|\kappa|=|\kappa\times\lambda|$.
\end{proof}

\begin{example}
We have $|\omega+1|=\omega=|\omega|$ even though $\omega \neq \omega +1$.
\end{example}
\begin{definition}
We write $\kappa^+$ for the least cardinal greater than $\kappa$.

We now use Transfinite Recursion to define the following:

    $\omega_0:=\omega$\\
    $\omega_{\alpha+1}:=(\omega_{\alpha})^+$\\
    If $\gamma$ is a limit ordinal  $\omega_{\gamma}:=sup\{\omega_{\alpha}: \alpha<\gamma\}$.
\end{definition}

\begin{lemma}\label{leminj}
    If $A\lesssim B$, then $|A|\leq|B|$.
\end{lemma}

\begin{theorem}\label{thmcantor} Cantor's Theorem.
    For any well-ordered set $x$ we have $|x|<|P(x)|$.
\end{theorem}

\begin{proof}
    There is a 1-1 map $\psi: x\rightarrow P(x)$ given by $a\rightarrow \{a\}$, hence $|x|\leq|P(x)|$ by Lemma \ref{leminj}. Suppose there is a bijection $\phi: x\rightarrow P(x)$. Define $S=\{a\in x: a\notin \phi(a)\}$. Then $S\in P(x)$ so there is an $s$ s.t. $\phi(s)=S$. If $s\in S$, then $s\notin S$ and if $s\in S$, then $s\notin S$. We reach a contradiction either way so there is no such bijection thus $|x|\neq|P(x)|$.
\end{proof}
By the Power Set Axiom and Comprehension, we can define:
\begin{definition}
    For any sets $A,B$ we define $^B A:=\{f\in P(B\times A): f \text{ is a function} \wedge dom(f)=B \wedge ran(f)\subseteq A\}$.

\end{definition}

\begin{definition}\label{defcardexp}
   Cardinal exponentiation is given by $\kappa^{\lambda}:=|^{\lambda}\kappa|$.
\end{definition}
The notation here is the same as for ordinal exponentiation \ref{defordexp} but these are different operations. In what follows we will always use this notation to refer to cardinal exponentiation.
\begin{example}
     The set $2^{\omega}$ is a cardinal.
\end{example}
\begin{lemma}\label{lem2kPk}
    For $\kappa \geq \omega$ we have $2^{\kappa}=|P(\kappa)|$.
\end{lemma}
This is my attempt at proving the previous lemma:
\begin{proof}
We define a bijection $\phi: P(\kappa)\rightarrow 
   ^{\omega}2$ by $\phi(S)=\{\langle x, 1\rangle :  x\in S \}\cup\{\langle y, 0 \rangle: y\in \kappa\backslash S\}$.
   
\end{proof}

From Lemma \ref{lem2kPk} and Theorem \ref{thmcantor} we have the following.
\begin{corollary}
    The strict inequality $\omega < 2^{\omega}$ holds.
\end{corollary}

\begin{definition}
    A limit ordinal $\kappa$  is \textit{regular} if every unbounded subset ${ C\subseteq \kappa }$ has cardinality $\kappa$. 
\end{definition}

\begin{lemma}
     
    A limit ordinal $\kappa$ is regular iff the union of $<\kappa$ subsets of $\kappa$ of cardinality $<\kappa$ has cardinality $<\kappa$.
\end{lemma}
\begin{corollary}
    The ordinal $\omega$ is regular.
\end{corollary}
\begin{proof}
    Any finite union of finite sets is finite.
\end{proof}

\section{Martin's Axiom}\label{sectionmartin}
We have shown that $\omega\neq 2^{\omega}$, now we are interested in whether $\omega^+=2^{\omega}$, that is, whether there are $\kappa$ s.t. $\omega<\kappa<2^{\omega}$.  In this section, we introduce \textit{Martin's Axiom} and then show that if we accept it, then the infinite cardinals $\kappa<2^{\omega}$ "behave" like $\omega$ in the sense that the answers to Questions 1-4 in this section are "yes" for $\omega \leq \kappa<2^{\omega}$. 

\begin{definition}
    The \textit{Continuum Hypothesis} (CH) is the statement $\omega_1=2^{\omega}$.
\end{definition}
CH is independent of ZFC so we could include it or its negation as an axiom of ZFC. If we reject CH, then there are infinite cardinals $\kappa$ between $\omega$ and $2^{\omega}$ and so, we are interested in the following questions:

Question 1: Is $2^{\omega}=2^{\kappa}$?\\
Question 2: Does every a.d. family of size$\kappa$ fail to be maximal?\\
Question 3: Does the union of $\kappa$ subsets of $\mathbb{R}$, each of Lebesgue measure 0, have measure 0?\\
Question 4: Does the union of $\kappa$ first-category subsets of $\mathbb{R}$ have first category?

We only discuss Questions 1 and 2.
\begin{definition}
    Two elements $x,y\subseteq \kappa$ for $\kappa$ an infinite cardinal are almost disjoint (a.d.) if $|x\cap y|<\kappa$. A family $A\subseteq P(\kappa)$ is a.d. if $\forall x\in A(|x|=\kappa)$ and any two $x,y\in A$ are a.d. An a.d. family $A$ is maximal if no a.d. family $B$ properly contains $A$.
\end{definition}
\begin{example}\label{evenodd}
    The sets $A,B$ of all even and odd natural numbers as subsets of $\omega$  are a.d. and $|A|=|B|=\omega$ hence $\{A,B\}\subset P(\omega)$ is an a.d. family.
    
\end{example}

\begin{theorem}\label{thmq2}
    The answer to Question 2 is "yes" for $\kappa$ regular, in particular, for $\kappa=\omega$.
\end{theorem}

\begin{proof}
    Let $A=\{A_{\xi}:\xi<\kappa\}$ be an a.d. family with $|A|=\kappa$. Let $B_{\xi}= A_{\xi}\backslash \bigcup_{\eta<\xi} A_{\eta}$ for each $\xi<\kappa$. Since the $A_{\xi}$ are a.d. $A_{\eta}\cap A_{\xi}<\kappa$, then by regularity of $\kappa$ we have $\bigcup_{\eta<\xi} (A_{\eta}\cap A_{\xi})<\kappa$. Thus $|B_{\xi}|= |A_{\xi}\backslash \bigcup_{\eta<\xi} A_{\eta}|=|A_{\xi}\backslash \bigcup_{\eta<\xi} (A_{\eta}\cap A_{\xi}|)\neq 0$ since $|A_{\xi}|=\kappa$. Thus $B_{\xi}\neq 0$ and we choose $\beta_{\xi}\in B_{\xi}$ for each $\xi< \kappa$. The $B_{\xi}$ are disjoint by construction so the $\beta_{\xi}$ are distinct. Let $D=\{\beta_{\xi}:\xi<\kappa\}$. Then if $\beta_{\xi}\in A_{\eta}$ we must have $\eta\geq\xi$. Thus $D\cap A_{\eta}\subseteq \{\beta_{\xi}:\xi<\eta\}$ which has cardinality $<\kappa$, hence $D, A_{\eta}$ are a.d. for each $\eta$, hence $A$ is properly contained in the a.d. family $A\cup \{D\}$.
\end{proof}

\begin{example}

Consider $n_0, n_1, n_2,...$ form Example \ref{exn0}

Again we use Transfinite Recursion to define the following:

$N^0:=\{n_0,n_1,...\}=\{0,1,3,6,10,...\}$\\
$N^i:=N^{i-1}+1$ where "+1" means that we add 1 to each coordinate.

For instance, $N^1=\{1,2,4,7,11, ...\}$ and $N^2=\{2,3,5,8,12, ...\}$.




It can be checked that each pair $N^i,N^j$ share at most $i-j$ elements for $i>j$. Then $N=\{N^{i}:i\in \omega\}$ is an a.d. family of subsets of $\omega$ with size $\omega$.

Let $N_i\backslash K_i$ be the set $N_i$ where we remove the set of elements $K_i$ which are contained in $N_{j}$ for some $j<i$. By the Axiom of Choice, we can make a collection $S$ with one element from each $N_i\backslash K_i$ and this will have size $\omega$ since  $N_i\backslash K_i\neq\emptyset$ for all $i$ (because $K_i$ is finite while $N_i$ is infinite). Then $S$ will be a.d. from each $N^i$ since each $N_i$ shares at most $i+1$ elements with $S$, thus the a.d. family $N$ is not maximal.

\end{example}

\begin{remark}
The condition of the theorem that the subsets have size $\kappa$ cannot be relaxed.
\end{remark}
This is my attempt at proving the previous remark.
\begin{proof}
    
    We give a counterexample. The set $S$ of all finite subsets of $\omega$ together with $\omega$ itself has size ${\omega}$ by Example \ref{exfinsubsets}. It is a.d. and is maximal since any subset $x$ of $\omega$ not contained in $S$ has infinite cardinality so $|x\cap \omega|=\omega$; thus $S$ is maximal.
    \end{proof}

  \begin{remark}
The condition of the theorem that the a.d. family must have size $\omega$ cannot be relaxed either.
\end{remark}
This is my attempt at proving the previous remark.

\begin{proof}
    Clearly $\{\omega\}$ is an a.d. family, and is maximal since any a.d. family containing $\omega$ can't contain any other infinite subset of $\omega$ which is a.d. from $\omega\in\{\omega\}$.
\end{proof}
\begin{remark}
    The converse of Theorem \ref{thmq2} doesn't hold.
\end{remark}
  
\begin{proof}
    By Zorn's Lemma, there is an a.d. maximal family $F$ of size $\geq\kappa^+$ (since any a.d. family of size $\kappa$ is properly contained in some other a.d. family), hence we can remove one element from $F$ so that it still has size $\geq\kappa^+$ and is still a.d. but is not maximal. 
\end{proof}
In order to state Martin's Axiom we give the following four definitions.

Throughout $P$ denotes a partial order.  The elements of P are called conditions and we say $p$ extends $q$ is $p\leq q$.

\begin{definition}
   Two elements $p,q\in P$ are \textit{compatible} if $\exists d\in P( d \leq p \wedge d\leq q)$. If $p,q$ are not compatible we say they are incompatible and denote it by $p\perp q$.
   
   A subset $S\subseteq P$ is an \textit{antichain} in $P$ if any two distinct elements in $P$ are incompatible.
   
   We say $P$ has the \textit{countable chain condition} (c.c.c.) if there is no uncountable antichain in $P$.
\end{definition}
\begin{example}
    In a Boolean algebra $\mathcal{B}$ two elements $a,b\in\mathcal{B}\backslash \{0\}$ are incompatible iff $a\wedge b=0$. 
\end{example}
\begin{example}\label{exK}
    Let $\mathcal{K}=\{p:p\text{ is a function}\wedge p\subset \omega \times 2\wedge |p|<\omega\}$. Let $p\leq q$ if $q\subseteq p$. Then $\langle \mathcal{K}, \leq \rangle$ is a partial order. Throughout we denote this partial order $\langle \mathcal{K}, \leq \rangle$.
    
    Then two elements $p,q\in \mathcal{K}$ are compatible iff $\forall x\in dom(p)\cap dom(q)(p(x)=q(x))$ since $p,q$ are both extensions of $p$ restricted to $dom(p)\cap dom(q)$. If $p,q$ are compatible, then $p\cup q$ is an extension of both. 
\end{example}

\begin{example}

Consider the partial order $P$ from, \ref{exfam}. Then two people $x,y$ are compatible iff they share a descendant, so the collection of all childless people is an anti-chain in $P$. Since there is a finite number of people, this has the c.c.c.
\end{example}

%
%
%
\end{example}
\begin{example}
    Let $P=\{(a,b): a<b, a,b\in \mathbb{R}\}$ and let $(a,b)\leq (c,d)$ if $[a,b]\subseteq (c,d)$. Let $S\subset P$ be an antichain, then for each interval $(a,b)\in S$ we can find a rational number $q\in (a,b)$ since $a\neq b$ and $\mathbb{Q}$ is dense (in the sense of topology rather than Definition \ref{defdense}) in $\mathbb{R}$. Thus $S$ is countable, since $\mathbb{Q}$ is, thus $P$ has the c.c.c. 
\end{example}    

More generally:
\begin{lemma}
    If $X$ is a separable topological space, then there is no uncountable family of pairwise disjoint open subsets of $X$.

    Where $X$ being \textit{separable} means that $X$ contains a countable, dense subset.
\end{lemma}
\begin{proof}
    Let $D$ be a countable and dense subset of $X$ (which exists since $X$ is separable). Let $\{U_{\alpha}\}$ be  a family of disjoint open subsets of $X$. Then for each $U_{\alpha}$ there is a $d\in D$ s.t. $d\in U_{\alpha}$. Hence since $D$ is countable, $\{U_{\alpha}\}$ is countable.
\end{proof}
\begin{definition}\label{defdense}
    A subset $A\subseteq P$ is dense in the partial order $P$ if $\forall p\in P (\exists d\in D (d\leq p))$.
\end{definition}

In particular, $P$ is dense in $P$ but we are interested in smaller dense sets. 
\begin{example}
    The subsets $\mathbb{Z},\mathbb{Q}\subseteq \mathbb{R}$ are dense in $\mathbb{R}$ with their usual orderings.
\end{example}
\begin{example}
    If we reverse the order in the partial order from Example \ref{exfam}, then the collection of all females is dense since everyone has a female ancestor.  
\end{example}
\begin{example}\label{exEh}
Consider the partial order $\mathcal{K}$ from Example \ref{exK}. Let $E_h:=\{f\in \mathcal{K}: \exists x\in dom(f): f(x)\neq h(x)\}$ for $h\in ^{\omega}2$. Let $f\in \mathcal{K}$. Since $dom(h)$ is infinite while $dom(f)$ is finite $\exists x\in dom(h)\backslash dom(f)$. Then we can extend $f$ to some $g$ with $dom(g)=dom(f)\cup\{x\}$ and define $g(x)\neq h(x)$. Then $g$ has finite domain and differs from $h$ so $g\in E_h$ and $g \leq f$. Thus $E_h$ is dense in $\mathcal{K}$. 

Let $D_n=\{f\in \mathcal{K}: n\in dom(f)\}$. Then $D_n$ is also dense since for any $f\in \mathcal{K}$ with $n\notin dom(f)$ we can extend $f$ to a function $g$ with $n\in dom(g)$.
 \end{example}

\begin{definition}
   A \textit{filter} $G$ on $P$ is a subset of $P$ s.t.:

1) $\forall p,q\in G (\exists r\in G(r \leq p \wedge r\leq q))$.\\
   2) $\forall q\in G (p \geq q \Rightarrow p\in G)$.
\end{definition}
\begin{example}
     The set of all ancestors of some (any) person is a filter in the partial order from \ref{exfam}.
\end{example}
\begin{definition}
    MA($\kappa$) is the statement: If $\langle P, \leq \rangle$ is a partial order of size $\kappa$ and $\mathcal{D}$ is a collection of $\leq \kappa$ dense subsets of $P$. Then there is a filter $G$ s.t. $\forall D\in \mathcal{D} (G\cap D\neq 0)$.

\end{definition}
\begin{axiom}
    Martin's Axiom. $MA(\kappa)$ for all $\kappa<2^{\omega}$.
\end{axiom}
\begin{example}
Any dense set in the partial order from Example \ref{exfam} has to include all childless people hence the filter consisting of all ancestors of some (any) childless person intersects any number of dense sets.



\end{example}

\begin{lemma}
    $MA(2^w)$ is not true.
\end{lemma}
\begin{proof}
    We give a counterexample:




    Consider $\mathcal{K}$ from Example \ref{exK}.
    

    Let $G$ be a filter on $P$. Then, since any two elements in $G$ are compatible we can define $f_G:= \bigcup G$ as the common extension of all functions in $G$.

    


Define $\mathcal{D}=\{\{D_n:n\in \omega\}\cup\{E_h:h\in ^{\omega}2\}\}$ where $E_h$ and $D_n$ are the dense sets from Example~\ref{exEh}. Then $|\mathcal{D}|=\omega+2^{\omega}=2^{\omega}$ by Theorem~\ref{thmmaxkl}.

    Now, if $\forall D\in \mathcal{D}(G\cap D\neq 0)$, then for every $n\in \omega$ there is a function $g\in G$ with $n\in dom(g)$ and since $f_G$ extends all of them $dom(f_G)=\omega$. Also, for every function $h\in ^{\omega}2$ there is a function $g\in G$ distinct from $h|_{dom(g)}$ and since $f_G$ extends all $g\in G$ we have that $f_G$ is different from $h$. Thus, $f_G\in ^{\omega}2$ and is distinct from all $h\in ^{\omega}2$, a contradiction. 
\end{proof}
\begin{lemma}
    $MA(\omega)$ is true.
\end{lemma}
\begin{proof}
    Let $\langle P, \leq \rangle$ be a partial order of size $\omega$ and let $D=\{D_n:n\in \omega\}$ be a family of dense subsets. Fix $p_0\in P$ and for each $n$ choose $p_n\in D_n$ s.t. $p_n\leq p_{n-1}$, which is possible since $D_n$ is dense. Define $G:=\{q\in P: \exists n(p_n\leq q)\}$, then $G$ is a filter and $\forall D_n\in D(G\cap D_n\neq0)$.
\end{proof}
Now we answer Questions 1 and 2 using the following partial order:

\begin{definition}
    For $A\subseteq P(\omega)$ the \textit{almost disjoint sets partial order} $P_A$ is $$\{\langle s, F\rangle: s\subset \omega \wedge |s|<\omega\wedge F\subseteq A \wedge |F|<\omega\}$$ where $\langle s', F'\rangle \leq \langle s, F\rangle$ if $s\subseteq s'\wedge F\subseteq F'\wedge \forall x \in F(x\cap s'\subseteq s)$.
\end{definition}

\begin{lemma}\label{lemcomp}
    Two elements $\langle s_1,F_1\rangle$ and $\langle s_2,F_2\rangle$ in $P_A$ are compatible iff
    $\forall x\in F_1( x\cap s_2\subseteq  s_1)\wedge 
    \forall x\in F_2(x\cap s_1\subseteq  s_2)$.
\end{lemma}
\begin{proof}
    Two elements $\langle s_1,F_1\rangle$ and $\langle s_2,F_2\rangle$ in $P_A$ are compatible \\$\iff \langle s_1\cup s_2,F_1\cup F_2\rangle$ extends $\langle s_1,F_1\rangle$ and $\langle s_2,F_2\rangle$  \\
    $\iff \forall x\in F_1(x\cap (s_1\cup s_2) \subseteq s_1) \wedge 
    \forall x\in F_2(x\cap (s_1\cup s_2) \subseteq s_2)$\\
    $\iff \forall x\in F_1(\forall n\in x\backslash s_1(n\notin s_1\cup s_2)\wedge 
    \forall x\in F_2(\forall n\in x\backslash s_2(n\notin s_1\cup s_2)$\\
    $\iff\forall x\in F_1(\forall n\in x\backslash s_1(n\notin  s_2)\wedge 
    \forall x\in F_2(\forall n\in x\backslash s_2(n\notin s_1)$\\
    $\iff\forall x\in F_1( x\cap s_2\subseteq  s_1)\wedge 
    \forall x\in F_2(x\cap s_1\subseteq  s_2)$.
\end{proof}
\begin{theorem}\label{thmUF}
    Assume MA(k). If $A,C\subseteq P(\omega)$ with $|A| \leq k, |B|\leq k$ and for all $y\in C$ and all finite $F$ $|y\backslash \bigcup F|=\omega$, then there exists $d\subseteq \omega$ s.t. $\forall x\in A(|d\cap x|<\omega)$ and $\forall x\in C(|d\cap x|=\omega)$
\end{theorem}
\begin{proof}
    For $y\in C$ and $n\in \omega$ let $E^y_n=\{\langle s,F\rangle\in P_A: s\cap y \not\subset n\}$. We first show $E^y_n$ is dense in $P_A$: Let $\langle s,F\rangle\in P_A$, then for $m\in y\backslash \bigcup F$ and $m>n$ we have $\langle s\cup {m},F\rangle$ extends $\langle s,F\rangle$ since $\{m\}\notin \bigcup F$ hence for $x\in F$ $x\cap s\cup\{m\}=x\cap s\subseteq s$. Also, $s\cup \{m\}\cap y \not\subseteq n$ since $m>n$, thus $\langle s\cup {m},F\rangle \in E^y_n$.

    For $x\in A$ let $D_x=\{\langle s,F\rangle: x\in F\}$. We show $D_x$ is dense: For $\langle s,F\rangle\in P_A$ we have $\langle s,F\cup\{x\}\rangle\in D_x$ and clearly extends $\langle s,F\rangle$.

    Let $D=\{\{D_x:x\in A\}\cup \{E^y_n:n\in \omega\wedge y\in C\}\}$. Then by MA(k), there is a filter $G$ s.t. $\forall D\in D(G\cap D\neq 0)$.

    Let $d_G=\bigcup\{s:\exists F(\langle s, F\rangle \in G)\}$. Then if $\langle s, F\rangle,\langle s', F'\rangle \in G$, then they are compatible hence $\forall x\in F(x\cap s'\subseteq s)$ by Lemma~\ref{lemcomp}. Thus, (since $\langle s', F'\rangle$ was arbitrary) $\forall x\in F(x\cap d_G\subseteq s)$. 

    Since $D_x\cap G\neq 0$ for $x\in A$ then $x\in F$ where $\langle s, F \rangle\in G$. Hence  $x\cap d_G\subseteq s$ and so $\forall x\in A(|x\cap d_G|<\omega)$.

    Also, for $y\in C$ we have $\forall n(y\cap d_G\not\subset n)$ hence $|y\cap d_G|=\omega$. So, $d=d_G$ works.
    \end{proof}
\begin{corollary}
Assume $MA(\kappa)$. Then the answer to Question 2 is "yes" for $\omega \leq \kappa < 2^{\omega}$.

\end{corollary}
\begin{proof}
Let $A\subset P(\omega)$ be an a.d. family with $|A|=\kappa$. We show $|\omega\backslash\bigcup F|=\omega$ for all finite $F\subset A$. Suppose $|\omega\backslash\bigcup F|<\omega$, then for $x\in A\backslash \bigcup F$ we have $|x|=\omega$ (because $A$ is an a.d. family) and $|x\cap \bigcup F|=\omega $ since $|x\cap \omega|=\omega $ and $|x\cap (\omega\backslash \bigcup F)|<\omega $. Since $F$ is finite there is some $y\in F$ s.t. $|x\cap y|=\omega $, contradicting $A$ being a.d. 

  Hence Theorem \ref{thmUF} applies with $C=\{\omega\}$ so there is a $d$ a.d. from all elements of $A$, so $A$ is not maximal.
\end{proof}

\begin{theorem}
    Assume $MA(k)$. Then the answer to Question 1 is "yes" for $\omega \leq \kappa <2^{\omega}$.
\end{theorem}

\begin{proof}
    For  $A\subset B\subseteq P(\omega)$ with $|B| = k$ where $|\omega\leq k<2^{\omega}|$ and $B$ is a.d. we have that Theorem~\ref{thmUF} applies with $A=A$ and $C=B\backslash A$. Hence there is a $d\subseteq \omega$ s.t. $\forall x\in A(|d\cap x|<\omega)$ and $\forall x\in B\backslash A(|d\cap x|=\omega)$ (*).

    Define $\phi: P(\omega)\rightarrow P(B)$ by $\phi(d)=\{x\in B: |x\cap d|<\omega\}$.

    By (*), we have that for any $A\in P(B)$ there is a $d\in P(\omega)$ s.t. $\phi(d)=A$, hence $\phi$ is onto. Thus $2^\omega=|P(\omega)|\geq |P(\kappa)|=2^{\kappa}$ by Lemma~\ref{lem2kPk}. Then the theorem follows from $2^\omega\leq 2^{\kappa}$.
\end{proof}

\section{Equivalents of MA}
In this last section we give equivalent forms of Martin's axiom, one of which is surprisingly a purely topological statement.

\begin{definition}
    Let $X$ be a topological space. Let $int(b)$ and $cl(b)$ be the interior and closure of $b$, respectively. 
    
    A subset $b\subseteq X$ is \textit{regular} (in the topological sense) if $b=int(cl(b))$. The regular open algebra $ro(X)$ of $X$ consists of the regular open subsets $b\subseteq X$ with operations given by $b\wedge c=b\cap c$, $v\vee c= int(cl(b \cup c)$ and the complement  of $b$ is $b'=int(X\backslash b)$.
\end{definition}
\begin{definition}
    A Boolean algebra is \textit{complete} if every subset has a supremum and an infimum. 
\end{definition}

\begin{lemma}\label{lembolalg}
    Let $P$ be a partial order. There exists a complete Boolean algebra $\mathcal{B}$ called the completition of $P$ and a map $i: P\rightarrow \mathcal{B}\backslash\{0\}$ s.t. :

    (1) $i'' P$ is dense in $\mathcal{B}\backslash\{0\}$.\\
    (2) $\forall p,q\in P$  $(p \leq q\Rightarrow i(p)\leq i(q))$.\\
    (3) $\forall p,q\in P$ $(p \perp q\iff i(p)\wedge i(q)=0)$.
    
\end{lemma}

\begin{proof}

Let $N_p=\{q\in P: q\leq p\}$ and define a topology on $P$ with base $\{N_p:p\in P\}$. Then define $\mathcal{B}=ro(P)$ and $i(p)=int(cl(N_p))$. If $q\in N_p$, then $1\leq p$ so $N_q\subseteq N_p$. We check these satisfy (1)-(3):

Let $b\neq 0\in ro(X)$. If $p\in b$, then $N_p\subseteq b$ so $i(p)=int(cl(N_p))\subseteq int(cl(b))=b$. Hence (1).

Let $p,q\in P$ s.t. $p\leq q$, then $i(p)=N_p\subseteq N_q=i(q)$. Hence (2).

If $p,q$ are compatible, then there is some $r\leq p,q$. Then $i(r)\leq i(p),i(q)$ thus $i(p)\wedge i(q) \neq 0$. Conversely, if $p,q$ are not compatible, then $N_p\cap N_q=0$, so $cl(N_p)\cap N_q=0$ since $N_q$ is open, so $i(p)\cap N_q=int (cl(N_p))\cap N_q=0$ since $int (cl(N_p))\subseteq cl(N_p))$. Then $i(p)\cap cl(N_q)=0$ since $i(p)$ is open, so $i(p)\cap i(q)=i(p)\cap int (cl(N_q))=0$. Hence (3).
\end{proof}

This gives us a way of associating a Boolean algebra to a partial order.
\begin{example}

    If any two elements in $P$ are compatible the completition $\mathcal{B}$ of $P$ is the 2-elements Boolean algebra: 

This is my attempt at showing this.

    Let $a\in i'' P$. Then $a'\notin i''P$ since $a \wedge a'=0$, so if they both had preimages in $i$, then they wouldn't be compatible by (3).

    Then either $a'=0$ or $\exists b\in i''P$ s.t. $b<a'$ by density of $i''P$ in $\mathcal{B}\backslash\{0\}$. If $\exists b\in i''P$ s.t. $b<a'$, then by compatibility of the preimages of $a$ and $b$, we have $a \wedge b=d\neq 0$. But then $d \leq a$ and $d\leq b <a'$ hence $a\wedge a' \geq d\neq 0$, a contradiction.

    Thus, $a'=0$ and then $a=1$, i.e., $\forall p\in P$ we have $i''(p)=1$. Then by density of $i''P$, we have that $\forall q\in \mathcal{B}\backslash \{0\}$ ($q\geq 1$), i.e., $\mathcal{B}$ is the Boolean algebra with two elements.
\end{example}

\begin{definition}
    An \textit{ultrafilter} of a partial order $P$ is a filter on $P$ which is not properly contained in any other filter on $P$.
\end{definition}

\begin{theorem}\label{equivMA}
    For any $\kappa\geq \omega$ TFAE:

    1) $MA(\kappa)$.\\
    2) $MA(\kappa)$ restricted to partial orders of cardinality $\leq \kappa$.\\
    3) $MA(\kappa)$ restricted to complete Boolean algebras, i.e. restricted to partial orders of the form $\mathcal{B}\backslash \{0\}$ where $\mathcal{B}$ is a complete c.c.c. Boolean algebra\\
    4) For any compact c.c.c. Haussdorf space $X$ and $U_{\alpha}$ are dense open subsets of $X$ for $\alpha<\kappa$, then $\cap_{\alpha} U_{\alpha}\neq 0$.
\end{theorem}
\begin{proof}
We prove only (4)$\Rightarrow$(3), (1)$\Rightarrow$(4) and (3)$\Rightarrow$(2), then it remains to show (2)$\Rightarrow$(1).

(4)$\Rightarrow$(3):

    Let $\mathcal{B}$ be a c.c.c. Boolean algebra (we don't require $\mathcal{B}$ to be complete). Let $\mathcal{D}$ be a collection of $\leq k$ dense subsets of $\mathcal{B}$. 

Define the following:

$X=$ the Stone space on $\mathcal{B}$, consisting of the ultrafilters on $\mathcal{B}$,\\
 $N_b=\{G\in X: b\in G\}$ for each $b\in \mathcal{B}$ and \\
$W_D=\bigcup \{N_b:b\in D\}$ for each $D\in \mathcal{D}$.

Then $X$ is Hausdroff, the $N_b$ are basic open sets in $X$ and for any $N_c, N_b\in X$ we have $N_c\cap N_b=0\iff b\wedge c= 0$ thus $X$ has the c.c.c. since $\mathcal{B}$ does.

The $W_D$ are open and dense in the topological sense: Since $D$ is dense in the order sense, there is an extension $b$ of $c$ for any $c\in \mathcal{B}$. Hence $b\wedge c\neq 0$ thus $N_c\cap N_b\neq 0$ and then $N_c\cap W_D\neq 0$. 

Hence we can apply (d) to define a filter $G\in \bigcap \{W_D:D\in \mathcal{D}\}$. Then for each $D\in \mathcal{D}$ we have $G\in W_D$, so $\exists b(G\in N_b)$, so $\exists b(b\in G)$, so $G\cap D \neq 0$.

For (1)$\Rightarrow$(4) we need the following.
\begin{definition}
    A non-empty family $X$ of subsets of a set $P$ has the \textit{finite intersection property} (FIP) if any finite subcollection of $X$ has non-empty intersection.
\end{definition}

Define a partial order $\langle P, \leq \rangle$ by $P=\{p\subseteq X: p\neq 0\textrm{ is open in } X\}$ and $p\leq q \iff p\subseteq q$. Then $p,q$ are incompatible iff $p\cap q=0$ and since $X$ has the c.c.c so does $P$. Let $G$ be a filter on $P$, then $G$ has the FIP, hence so does $\{\bar{b}:b\in G\}$ since $b\subseteq \bar{b}$. A topological space is compact iff any collection of its closed sets having the FIP has non-empty intersection. Hence $\bigcap\{\bar{b}:b\in G\}\neq 0$. Then for each $\alpha$ the set $D_{\alpha}:=\{p\in P: \bar{p}\in U_\alpha\} $ is dense in $P$ since $U_{\alpha}$ is dense in $X$ and $X$ is regular. By $MA(\kappa)$, there is a filter $G$ s.t. $\forall \alpha<\kappa(G\cap D_{\alpha}\neq 0)$. Then $0\neq \bigcap\{\bar{b}:b\in G\}\subseteq \bigcap_{\alpha} U_{\alpha}\neq 0$.

(3)$\Rightarrow$(2): 

Let $P$ be a partial order with $|P|\leq \kappa$ and $\mathcal{D}$ a collection of dense subsets of $P$. Let $\mathcal{B}$ and $i$ satisfy the conditions of Lemma \ref{lembolalg}. 

For $b\in \mathcal{B}\backslash \{0\}$ there is $i(p)\leq b$ since $i''P$ is dense in $\mathcal{B}$ by  Lemma \ref{lembolalg} (1). Then there is a $d\in D$ s.t. $d\leq p$ since $D$ is dense in $P$. Hence $i(d)\leq i(p) \leq b$, so $i''D$ is dense in $\mathcal{B}$.

Hence we can apply $MA(k)$ restricted to Boolean algebras to get a filter $G$ intersecting $i''D$ for each $D\in \mathcal{D}$. Then for $H:=i^{-1}(G)$ we have $H\cap D\neq 0$. 

Now we check that $H$ is a filter.

If $ P\ni p \geq h\in H$, then $\mathcal{B}\ni i(p)\geq i(h)\in G$ by  Lemma \ref{lembolalg} (2) and since $G$ is a filter $i(p)\in G$, thus $p\in H$. 

Let $D_{pq}=\{r\in P: (r\leq p\wedge r\leq q)\vee r\perp p\vee r\perp q\}$ for all $p,q\in P$. Then each $D_{pq}$ is dense in $P$; Let $r_0\in P$. If $\exists r_1\leq r_0 $ s.t. $r_1\perp p\vee r_1\perp q$, then $r_1\in D_{pq}$. If not, then every extension of $r_0$ is compatible with $p$ and $q$, in particular, $r_0$ is. Then $\exists r_1 \leq r_0, p$, so $r_1$ is compatible with $q$, so we have $\exists r_2 \leq r_1, q$. Then $r_2 \leq p,q, r_1$ so $r_2\in D_{pq}$ extends $r_0$.

Since $|P|\leq \kappa$ we have $|\mathcal{D}\cup\{D_{pq}:p,q\in P\}|\leq\kappa$ so we assume $D_{pq}\in \mathcal{D}$. Thus $\exists r\in H\cap D_{pq}$ and since elements of $H$ are pairwise compatible by Lemma \ref{lembolalg} (3), we have that $r\leq p,q$ and $r\in H$. Thus $H$ is a filter.   
\end{proof}
\newpage
\printbibliography

@book{Kunen, place={}, title={Set theory: An Introduction to Independence Proofs}, publisher={}, author={Kunen, Kenneth}, year={2006} }

@misc{russells,
author       =	{Irvine, Andrew David and Deutsch, Harry},
title        =	{{Russell's Paradox}},
booktitle    =	{The {Stanford} Encyclopedia of Philosophy},
editor       =	{Edward N. Zalta},
howpublished =	{\url{https://plato.stanford.edu/archives/spr2021/entries/russell-paradox/}},
year         =	{2021},
edition      =	{{S}pring 2021},
publisher    =	{Metaphysics Research Lab, Stanford University}
}

@article{martin,
	author = {Martin, D.A. and Solovay, R.M.},
	doi = {10.1016/0003-4843(70)90009-4},
	journal = {Annals of Mathematical Logic},
	%number = {1},
	title = {{Internal Cohen Extensions}},
	volume = {2},
	year = {1970},
}

@article{zfc, 
   author = {Zermelo, Ernst}, 
   title = {{Investigations in the foundations of set theory I}}, 
   publisher = {B. G. Teubner Verlag}, 
   journal = {Mathematische Annalen}, 
   volume = {65}, 
   number = {2}, 
   year = {1908}, 
   address = {Leipzig},  
 }

@article{Cantor,
author = {Cantor, Georg},
journal = {Jahresbericht der Deutschen Mathematiker-Vereinigung},
pages = {72-78},
title = {On an Elementary Question on Set Theory},
%url = {http://eudml.org/doc/144383},
volume = {1},
year = {1890/91},
}

@article{cohen,
% ISSN = {00278424},
% URL = {http://www.jstor.org/stable/71858},
 author = {Paul J. Cohen},
 journal = {Proceedings of the National Academy of Sciences of the United States of America},
 number = {6},
% pages = {1143--1148},
 publisher = {National Academy of Sciences},
 title = {The Independence of the Continuum Hypothesis},
% urldate = {2023-01-17},
 volume = {50},
 year = {1963}
}


\end{document}